\documentclass{amsart}

\usepackage{amsmath,amsfonts,amssymb}

\usepackage{tikz}
\usepackage{multirow}
\usepackage{graphicx}
\usepackage{verbatim}
\usepackage{url}
\usepackage{hyperref}

\usepackage{float}
\floatstyle{plaintop}
\restylefloat{table}

\theoremstyle{plain}
\newtheorem{theorem}{Theorem}[section]

\newtheorem{proposition}[theorem]{Proposition}

\theoremstyle{definition}

\theoremstyle{remark}

\newcommand{\C}{\mathbb{C}}

\begin{document}

\newpage
\title{Explicit formulae for Chern-Simons invariants of the hyperbolic orbifolds of the knot with Conway's notation $C(2n, 3)$}
\author{Ji-Young Ham, Joongul Lee}
\address{Department of Science, Hongik University, 
94 Wausan-ro, Mapo-gu, Seoul,
 04066\\
   Department of Mathematical Sciences, Seoul National University, 
   1 Gwanak-ro, Gwanak-gu, Seoul 08826 \\
   Korea} 
\email{jiyoungham1@gmail.com.}

\address{Department of Mathematics Education, Hongik University, 
94 Wausan-ro, Mapo-gu, Seoul,
 04066\\
   Korea} 
\email{jglee@hongik.ac.kr}



\maketitle 

\markboth{ Ji-young Ham, Joongul Lee } 
{ Chern-Simons invariants of $C(2n,3)$ }
 
\begin{abstract}
 We calculate the Chern-Simons invariants of  the hyperbolic orbifolds of the knot with Conway's notation $C(2n, 3)$ using the Schl\"{a}fli formula for the generalized Chern-Simons function on the family of $C(2n,3)$ cone-manifold structures. We present the concrete and explicit formula of them. We apply the general instructions of Hilden, Lozano, and Montesinos-Amilibia and extend the Ham and Lee's methods.
As an application, we  calculate the Chern-Simons invariants of cyclic coverings of the hyperbolic $C(2n,3)$ orbifolds.
\end{abstract}
\maketitle

\section{Introduction}
Chern-Simons invariant~\cite{CS, Mey1} was defined to be a geometric invariant and became a topological invariant 
after the Mostow Rigidity Theorem~\cite{Mo1}.
Various methods of finding Chern-Simons invariant using ideal triangulations have been 
introduced~\cite{N1,N2,Z1,CMY1,CM1,CKK1} and
 implemented~\cite{SnapPy, Snap}. But, for orbifolds, to our knowledge, there does not exist a single convenient program which computes Chern-Simons invariant.

Instead of working on complicated combinatorics of 3-dimensional ideal tetrahedra to find the Chern-Simons invariants of the hyperbolic orbifolds of the knot with Conway's notation $C(2n, 3)$, we deal with simple one dimensional singular loci. 
Similar methods for volumes can be found in~\cite{HL1, HMP}. We use the Schl\"{a}fli formula for the generalized Chern-Simons function on the family of $C(2n,3)$ cone-manifold structures~\cite{HLM3}.  
In~\cite{HLM2} a method of calculating the Chern-Simons invariants of two-bridge knot orbifolds were introduced but without explicit formulae. In~\cite{HL}, the Chern-Simons invariants of the twist knot orbifolds are computed.
Similar approaches for $SU(2)$-connections can be found in~\cite{KK1} and for $\text{\textnormal{SL}}(2,C)$-connections in~\cite{KK2}.
For explanations of cone-manifolds, you can refer to~\cite{CHK,T1,K1,P2,HLM1,PW,HMP}.

\medskip 

The main purpose of the paper is to find the explicit and efficient formulae for  Chern-Simons invariants of the hyperbolic orbifolds of the knot with Conway's notation $C(2n, 3)$.  For a two-bridge hyperbolic link, there exists an angle $\alpha_0 \in [\frac{2\pi}{3},\pi)$ for each link $K$ such that the cone-manifold $K(\alpha)$ is hyperbolic for $\alpha \in (0, \alpha_0)$, Euclidean for $\alpha=\alpha_0$, and spherical for $\alpha \in (\alpha_0, \pi]$ \cite{P2,HLM1,K1,PW}. We will use the Chern-Simons invariant of the lens space $L(6n+1,4n+1)$ calculated in~\cite{HLM2}.
The following theorem gives the formulae  for $T_{2 n}$. Note that if $2 n$ of $T_{2n}$ is replaced by an odd integer, then $T_{2n}$ becomes a link with two components. Also, note that the Chern-Simons invariant of hyperbolic cone-manifolds of the knot with Conway's notation $C(-2n, -3)$ is the same as that of the knot with Conway's notation $C(2n,3)$ up to sign. For the Chern-Simons invariant formula, since the knot $T_{2n}$ has to be hyperbolic, we exclude the case when $n=0$.

\begin{theorem}\label{theorem:main}
Let $X_{2n}(\alpha)$, $0 \leq \alpha < \alpha_0$ be the hyperbolic cone-manifold with underlying space $S^3$ and with singular set $T_{2n}$ of cone-angle $\alpha$. Let $k$ be a positive integer such that $k$-fold cyclic covering of $X_{2n}(\frac{2 \pi}{k})$ is hyperbolic. Then the Chern-Simons invariant of $X_{2n}(\frac{2 \pi}{k})$ (mod $\frac{1}{k}$ if $k$ is even or mod $\frac{1}{2k}$ if $k$ is odd) is given by the following formula:

\begin{equation*}
\begin{split}
&\text{\textnormal{cs}} \left(X_{2n} \left(\frac{2 \pi}{k} \right)\right) 
 \equiv \frac{1}{2} \text{\textnormal{cs}}\left(L(6n+1,4n+1) \right) \\
&+\frac{1}{4 \pi^2}\int_{\frac{2 \pi}{k}}^{\alpha_0} Im \left(2*\log \left(-M^{-4n-2}\frac{M^{-2}+t}{ M^{2}+t}\right)\right) \: d\alpha \\
& +\frac{1}{4 \pi^2}\int_{\alpha_0}^{\pi}
 Im \left(\log \left(-M^{-4n-2}\frac{M^{-2}+t_1}{ M^{2}+t_1}\right)+\log \left(-M^{-4n-2}\frac{M^{-2}+t_2}{ M^{2}+t_2}\right)\right) \: d\alpha,
\end{split}
\end{equation*}

\noindent where 
for $M=e^{\frac{\alpha}{2}}$, $t$ $(Im(t) \leq 0)$, $t_1$, and $t_2$ are zeroes of Riley-Mednykh polynomial $P_{2n}=P_{2n}(t,M)$ which is
 given recursively by 

\medskip
\begin{equation*}
P_{2n} = \begin{cases}
 Q P_{2(n-1)} -M^8 P_{2(n-2)} \ \text{if $n>1$}, \\
 Q P_{2(n+1)}-M^8 P_{2(n+2)} \ \text{if $n<-1$},
\end{cases}
\end{equation*}
\medskip

\noindent with initial conditions
\begin{equation*}
\begin{split}
& P_{-2}  =M^2 t^2+\left(M^4- M^2+1\right) t+M^2,\\
& P_{0}  =M^6 \ \text{\textnormal{for}} \ n \leq 0 \qquad \text{\textnormal{and}} \qquad P_{0}  = M^8 \ \text{\textnormal{for}} \ n \geq 0,\\    
& P_{2}  =-M^4 t^3+\left(-2M^6+M^4-2 M^2\right) t^2+\left(-M^8+M^6-2 M^4+ M^2-1\right) t+M^4, \\
\end{split}
\end{equation*}
\noindent and $M=e^{\frac{\alpha}{2}}$ and
$$Q=-M^4 t^3+\left(-2M^6+2 M^4-2M^2\right) t^2+\left(-M^8+2M^6-3M^4+2M^2-1\right) t+2M^4,$$

\medskip

\noindent where $M=e^{\frac{\alpha}{2}}$ and $t_1$ and $t_2$ approach common $t$ as $\alpha$ decreases to $\alpha_0$ and they come from the components of $t$ and $\overline{t}$.
\end{theorem}

%%%%%%%%%%%%%%%%%%%%%%%%%%%%%%%%%%%%%%%%%%%%%%%%%%%%%%%%
%%%%%%%%%%%%%%%%%%%%%%%%%%%%%%%%%%%%%%%%%%%%%%%%%%%%%%%%

\section{Two bridge knots with  Conway's notation $C(2n, 3)$} \label{sec:C[2n,3]}

\begin{figure} 
\begin{center}
\resizebox{3.5cm}{!}{\includegraphics[angle=90]{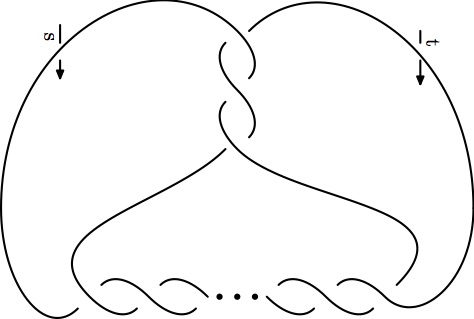}}
reflection
\reflectbox{\resizebox{3.5cm}{!}{\includegraphics[angle=90]{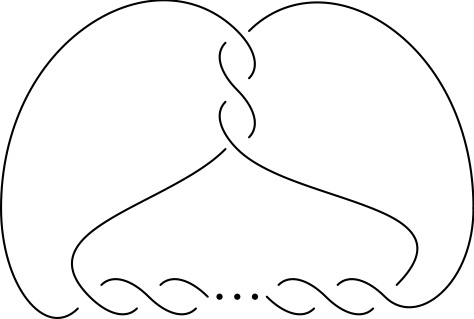}}}
\caption{A two bridge knot with Conway's notation C[2n,3] (left) and its mirror image C[-2n,-3](right)} \label{fig:C[2n,3]}
\end{center}
\end{figure}

\begin{figure}
\begin{center}
\resizebox{3.5cm}{!}{\includegraphics[angle=90]{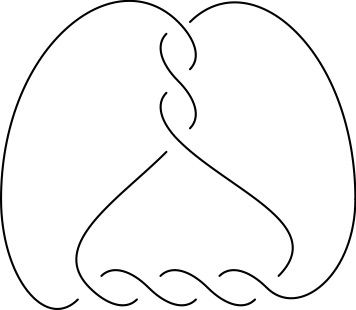}}
\caption{The knot $7_3$}\label{fig:knot}
\end{center}
\end{figure}

A knot $K$ is a two bridge knot with Conway's notation $C(2n,3)$ if $K$ has a regular two-dimensional projection of the form in Figure~\ref{fig:C[2n,3]}. For example, Figure~\ref{fig:knot} is knot $C(4,3)$.
$K$ has 3 left-handed horizontal crossings and $2 n$ right-handed vertical crossings.  
We will denote it by $T_{2 n}$. One can easily check that the slope of $T_{2n}$ is $3/(6n+1)$ which is equivalent to the knot with slope 
$(4n+1)/(6n+1)$~\cite{S1}. For example, Figure~\ref{fig:knot}  shows the regular projections of knot $7_3$ with slope $3/13$ which is equivalent to the knot with slope $9/13$ (right).

We will use the following fundamental group of the knot with Conway's notation $C(2n, 3)$~\cite{HL1,HS,R1}. The following theorem can also be obtained by reading off the fundamental group from the Schubert normal form of $T_{2n}$ with slope $\frac{4n+1}{6n+1}$~\cite{S1,R1}.

\begin{proposition}\label{theorem:fundamentalGroup}
$$\pi_1(X_{2n})=\left\langle s,t \ |\ swt^{-1}w^{-1}=1\right\rangle,$$
where $w=(ts^{-1}tst^{-1}s)^n$.
\end{proposition}
%%%%%%%%%%%%%%%%%%%%%%%%%%%%%%%%%%%%%%%%%%%%%%%%%%%%%%%%%%%%%
%%%%%%%%%%%%%%%%%%%%%%%%%%%%%%%%%%%%%%%%%%%%%%%%%%%%%%%%%%%%%
\section{The Riley-Mednykh polynomial}
Given  a set of generators, $\{ s,t \}$, of the fundamental group for 
$\pi_1 (X_{2n})$, we define a representation $\rho \ : \ \pi_1 (X_{2n}) \rightarrow \text{SL}(2, \C)$ by

\begin{center}
$$\begin{array}{ccccc}
\rho(s)=\left[\begin{array}{cc}
                       M &       1 \\
                        0      & M^{-1}  
                     \end{array} \right]                          
\text{,} \ \ \
\rho(t)=\left[\begin{array}{cc}
                   M &  0      \\
                   2-M^2-M^{-2}-t      & M^{-1} 
                 \end{array}  \right].
\end{array}$$
\end{center}

Then $\rho$ can be identified with the point $(M,t) \in \C^2$.  By~\cite{HL1}, when $M$ varies we have an algebraic set 
whose defining equation is the following Riley-Mednykh polynomial.

\begin{theorem} \label{thm:RMpolynomial}
$t$ is a root of the following Riley-Mednykh polynomial $P_{2n}=P_{2n}(t,M)$ which is given recursively by 

\medskip
\begin{equation*}
P_{2n} = \begin{cases}
 Q P_{2(n-1)} -M^8 P_{2(n-2)} \ \text{if $n>1$}, \\
 Q P_{2(n+1)}-M^8 P_{2(n+2)} \ \text{if $n<-1$},
\end{cases}
\end{equation*}
\medskip

\noindent with initial conditions
\begin{align*}
& P_{-2}   =M^2 t^2+\left(M^4- M^2+1\right) t+M^2,\\
& P_{0}  =M^6 \ \text{\textnormal{for}} \ n \leq 0 \qquad \text{\textnormal
{and}} \qquad P_{0}  =M^8 \ \text{\textnormal{for}} \ n \geq 0,\\    
& P_{2} (t,M)  =-M^4 t^3+\left(-2M^6+M^4-2 M^2\right) t^2+\left(-M^8+M^6-2 M^4+ M^2-1\right) t+M^4, \\
\end{align*}
\noindent and
$$Q=-M^4 t^3+\left(-2M^6+2 M^4-2M^2\right) t^2+\left(-M^8+2M^6-3M^4+2M^2-1\right) t+2M^4.$$
\end{theorem}

%%%%%%%%%%%%%%%%%%%%%%%%%%%%%%%%%%%%%%%%%%%%%%%%%%%%%%%%%%%%%
%%%%%%%%%%%%%%%%%%%%%%%%%%%%%%%%%%%%%%%%%%%%%%%%%%%%%%%%%%%%%
\subsection{Longitude}
\label{sec:longitude}
   Let $l = ww^{*}M^{-4n}$, where $w^{*}$ is the word obtained by reversing $w$. Let $L=\rho(l)_{11}$. Then $l$ is the longitude which is null-homologus in $X_{2n}$. And we have
   
\medskip
   
\begin{theorem}~\cite{HL1}
\label{thm:longitude}
\begin{align*}
L=-M^{-4n-2}\frac{M^{-2}+t}{ M^{2}+t}.
\end{align*}
\end{theorem}
%%%%%%%%%%%%%%%%%%%%%%%%%%%%%%%%%%%%%%%%%%%%%%%%%%%%%%%%
%%%%%%%%%%%%%%%%%%%%%%%%%%%%%%%%%%%%%%%%%%%%%%%%%%%%%%%%
\section{Schl\"{a}fli formula for the generalized Chern-Simons function}
\label{sec:CSfunction}

The general references for this section are~\cite{HLM3,HLM2,Y1,MeyRub1} and~\cite{HL}. 
We introduce the generalized Chern-Simons function on the family of $C(2n,3)$ cone-manifold structures. For the oriented knot $T_{2n}$, we orient a chosen meridian $s$ such that the orientation of $s$ followed by orientation of $T_{2n}$ coincides with orientation of $S^3$. Hence, we use the definition of Lens space in~\cite{HLM2} so that we can have the right orientation when the definition of Lens space is combined with the following frame field. On the Riemannian manifold $S^3-T_{2n}-s$ we choose a special frame field $\Gamma$. A \emph{special} frame field $\Gamma=(e_1,e_2,e_3)$ is an orthonomal frame field such that for each point $x$ near $T_{2n}$, $e_1(x)$ has the knot direction, $e_2(x)$ has the tangent direction of a meridian curve, and $e_3(x)$ has the knot to point direction.  A special frame field always exists by Proposition $3.1$ of~\cite{HLM3}. From $\Gamma$ we obtain an orthonomal frame field 
$\Gamma_{\alpha}$ on $X_{2n}(\alpha)-s$ by the Schmidt orthonormalization process with respect to the geometric structure of the cone manifold $X_{2n}(\alpha)$. Moreover it can be made special by deforming it in a neighborhood of the singular set and 
$s$ if necessary. $\Gamma^{\prime}$ is an extention of $\Gamma$ to $S^3-T_{2n}$. For each cone-manifold $X_{2n}(\alpha)$, we assign the real number:

\begin{equation*}
I\left(X_{2n}(\alpha)\right)=\frac{1}{2} \int_{\Gamma(S^3-T_{2n}-s)}Q-\frac{1}{4 \pi} \tau(s,\Gamma^{\prime})-\frac{1}{4 \pi} \left(\frac{\beta \alpha}{2 \pi}\right),
\end{equation*}

\noindent where $-2 \pi \leq \beta \leq 2 \pi$, $Q$ is the Chern-Simons form:

\begin{equation*}
Q=\frac{1}{4 \pi^2} \left(\theta_{12} \wedge \theta_{13} \wedge \theta_{23} + \theta_{12} \wedge \Omega_{12} + \theta_{13} \wedge \Omega_{13} + \theta_{23} \wedge \Omega_{23} \right),
\end{equation*}

\noindent and 

\begin{equation*}
\tau(s,\Gamma^{\prime})=-\int_{\Gamma^{\prime}(s)} \theta_{23},
\end{equation*}

\noindent where ($\theta_{ij}$) is the connection $1$-form, ($\Omega_{ij}$) is the curvature $2$-form of the Riemannian connection on $X_{2n}(\alpha)$ and the integral is over the orthonomalizations of the same frame field. When $\alpha = \frac{2 \pi}{k}$ for some positive integer, 
$I \left(X_{2n}\left(\frac{2 \pi}{k}\right)\right)$ (mod $\frac{1}{k}$ if $k$ is even or mod $\frac{1}{2k}$ if $k$ is odd) is independent of the frame field $\Gamma$ and of the representative in the equivalence class $\overline{\beta}$ and hence an invariant of the orbifold $X_{2n}\left(\frac{2 \pi}{k}\right)$. $I \left(X_{2n}\left(\frac{2 \pi}{k}\right)\right)$ (mod $\frac{1}{k}$ if $k$ is even or mod $\frac{1}{2k}$ if $k$ is odd) is called \emph{the Chern-Simons invariant of the orbifold} and is denoted by 
$\text{\textnormal{cs}} \left(X_{2n}\left(\frac{2 \pi}{k}\right) \right)$.

On the generalized Chern-Simons function on the family of $C(2n,3)$ cone-manifold structures we have the following Schl\"{a}fli formula.

\begin{theorem}(Theorem 1.2 of~\cite{HLM2})~\label{theorem:schlafli}
For a family of geometric cone-manifold structures, $X_{2n}(\alpha)$, and differentiable functions $\alpha(t)$ and $\beta(t)$ of $t$ we have
\begin{equation*}
dI \left(X_{2n}(\alpha)\right)=-\frac{1}{4 \pi^2} \beta d \alpha.
\end{equation*} 
\end{theorem}
%%%%%%%%%%%%%%%%%%%%%%%%%%%%%%%%%%%%%%%%%%%%%%%%%%%%%%%%
%%%%%%%%%%%%%%%%%%%%%%%%%%%%%%%%%%%%%%%%%%%%%%%%%%%%%%%%

\section{Proof of the theorem~\ref{theorem:main}} \label{sec:proof}

For $n \geq 1$ and $M=e^{i \frac{\alpha}{2}}$, $P_{2n}(t,M)$ have $3n$ component zeros, and for  $n \leq -1$, $-(3n+1)$ component zeros. The component which gives the maximal volume is the geometric component~\cite{D1,FK1,HL1} and in~\cite{HL1} it is identified. For each $T_{2n}$, there exists an angle $\alpha_0 \in [\frac{2\pi}{3},\pi)$ such that $T_{2n}$ is hyperbolic for $\alpha \in (0, \alpha_0)$, Euclidean for $\alpha=\alpha_0$, and spherical for $\alpha \in (\alpha_0, \pi]$ \cite{P2,HLM1,K1,PW}.
Denote by $D(X_{2n}(\alpha))$ be the set of  zeros of the discriminant of 
$P_{2n}(t,e^{i \frac{\alpha}{2}})$ over $t$. Then $\alpha_0$ will be one of $D(X_{2n}(\alpha))$. 

 On the geometric component we can calculate
 the Chern-Simons invariant of an orbifold 
$X_{2n}(\frac{2 \pi}{k})$ (mod $\frac{1}{k}$ if $k$ is even or mod $\frac{1}{2k}$ if $k$ is odd), where $k$ is a positive integer such that $k$-fold cyclic covering of $X_{2n}(\frac{2 \pi}{k})$ is hyperbolic:
\begin{align*}
& \text{\textnormal{cs}}\left(X_{2n} \left(\frac{2 \pi}{k} \right)\right) 
                        \equiv I \left(X_{2n} \left(\frac{2 \pi}{k} \right)\right) 
                      \ \ \ \ \ \ \ \ \ \ \ \  \left(\text{mod} \ \frac{1}{k}\right) \\
                        & \equiv I \left(X_{2n}( \pi) \right)
                          +\frac{1}{4 \pi^2}\int_{\frac{2 \pi}{k}}^{\pi} \beta \: d\alpha 
                      \ \ \ \ \ \ \ \ \ \ \ \ \  \ \ \ \ \ \ \  \left(\text{mod} \ \frac{1}{k}\right) \\
                        & \equiv \frac{1}{2} \text{\textnormal{cs}}\left(L(6n+1,4n+1) \right) 
+\frac{1}{4 \pi^2}\int_{\frac{2 \pi}{k}}^{\alpha_0} Im \left(2*\log \left(-M^{-4n-2}\frac{M^{-2}+t}{ M^{2}+t}\right)\right) \: d\alpha \\
& +\frac{1}{4 \pi^2}\int_{\alpha_0}^{\pi}
 Im \left(\log \left(-M^{-4n-2}\frac{M^{-2}+t_1}{ M^{2}+t_1}\right)+\log \left(-M^{-4n-2}\frac{M^{-2}+t_2}{ M^{2}+t_2}\right)\right) \: d\alpha \\
& \left( \text{mod} \ \frac{1}{k}\ \text{if $k$ is even or }  \text{mod} \ \frac{1}{2k}\ \text{if $k$ is odd} \right)
\end{align*}
where the second equivalence comes from Theorem~\ref{theorem:schlafli} and the third equivalence comes from the fact that $I \left(X_{2n}(\pi)\right) \equiv \frac{1}{2} \text{\textnormal{cs}}\left(L(6n+1,4n+1) \right)$  $\left(\text{mod }\frac{1}{2}\right)$, Theorem~\ref{thm:longitude}, and geometric interpretations of hyperbolic and spherical holonomy representations.

The following theorem gives the Chern-Simons invariant of the Lens space $L(6n+1,4n+1)$.

\begin{theorem}(Theorem 1.3 of~\cite{HLM2}) \label{theorem:Lens}
\begin{align*}
\text{\textnormal{cs}} \left(L \left(6n+1,4n+1\right)\right) \equiv \frac{4n+4}{12n+2} && (\text{mod}\ 1).
\end{align*}
\end{theorem}

%%%%%%%%%%%%%%%%%%%%%%%%%%%%%%%%%%%%%%%%%%%%%%%%%%%%%%%%%
%%%%%%%%%%%%%%%%%%%%%%%%%%%%%%%%%%%%%%%%%%%%%%%%%%%%%%%%%
\section{Chern-Simons invariants of the hyperbolic orbifolds of the knot with Conway's notation $C(2n, 3)$ and of its cyclic coverings}

The table~\ref{table1-1} (resp. the table~\ref{table1-2}) gives the approximate Chern-Simons invariant of the hyperbolic orbifold, 
$\text{\textnormal{cs}} \left(X_{2n} (\frac{2 \pi}{k})\right)$ for $n$ between $2$ and $9$ (resp. for $n$ between $-9$ and $-2$) and for $k$ between $3$ and $10$, and of its cyclic covering, $\text{\textnormal{cs}} \left(M_k (X_{2n})\right)$. We used Simpson's rule for the approximation with $10^4$ ($5 \times 10^3$ in Simpson's rule) intervals from $\frac{2 \pi}{k}$ to $\alpha_0$ and $10^4$ ($5 \times 10^3$ in Simpson's rule) intervals from $\alpha_0$ to $\pi$. 
The table~\ref{tab2} gives the approximate Chern-Simons invariant of 
$T_{2n}$ for each n between $-9$ and $9$ except the unknot and the amphicheiral knot.  We again used Simpson's rule for the approximation with $10^4$ ($5 \times 10^3$ in Simpson's rule) intervals from $0$ to $\alpha_0$ and $10^4$ ($5 \times 10^3$ in Simpson's rule) intervals from $\alpha_0$ to $\pi$. 
We used Mathematica for the calculations. We record here that our data in table~\ref{tab2} and those obtained from  SnapPy match up up to six decimal points.

\begin{table}
\begin{tabular}{ll}
\begin{tabular}{|c|c|c|}
\hline
 $k$ & $\text{\textnormal{cs}} \left(X_2(\frac{2 \pi}{k})\right)$ & $\text{\textnormal{cs}} \left(M_k( X_2)\right)$ \\
\hline
  3 & 0.0200137 \text{} & 0.0600411 \text{} \\
 4 & 0.186810 & 0.747239 \\
 5 & 0.00166425 \text{} & 0.00832123 \text{} \\
 6 & 0.0504594 & 0.302756 \\
 7 & 0.0163411 \text{} & 0.114387 \text{} \\
 8 & 0.116987 & 0.935894 \\
 9 & 0.0292866 \text{} & 0.263580 \text{} \\
 10 & 0.0595395 & 0.595395 \\
\hline
\end{tabular}
&
\begin{tabular}{|c|c|c|}
\hline
 $k$ &  $\text{\textnormal{cs}} \left(X_4(\frac{2 \pi}{k})\right)$ & $\text{\textnormal{cs}} \left(M_k (X_4)\right)$ \\
\hline
 3 & 0.163905 \text{} & 0.491714 \text{} \\
 4 & 0.207480 & 0.829920 \\
 5 & 0.0602662 \text{} & 0.301331 \text{} \\
 6 & 0.140577 & 0.843464 \\
 7 & 0.0610011 \text{} & 0.427008 \text{} \\
 8 & 0.00457501 & 0.0366000 \\
 9 & 0.0181733 \text{} & 0.163560 \text{} \\
 10 & 0.0302655 & 0.302655 \\
\hline
\end{tabular}
\end{tabular}

\bigskip

\begin{tabular}{ll}
\begin{tabular}{|c|c|c|}
\hline
 $k$ &   $\text{\textnormal{cs}} \left(X_6(\frac{2 \pi}{k})\right)$ & $\text{\textnormal{cs}} \left(M_k (X_6)\right)$ \\
\hline
 3 & 0.0117308 \text{} & 0.0351925 \text{} \\
 4 & 0.0254160 & 0.101664 \\
 5 & 0.0770172 \text{} & 0.385086 \text{} \\
 6 & 0.130155 & 0.780930 \\
 7 & 0.0343996 \text{} & 0.240797 \text{} \\
 8 & 0.0925471 & 0.740377 \\
 9 & 0.0295838 \text{} & 0.266254 \text{} \\
 10 & 0.0810442 & 0.810442 \\
\hline
\end{tabular}
&
\begin{tabular}{|c|c|c|}
\hline
 $k$ &  $\text{\textnormal{cs}} \left(X_{8} (\frac{2 \pi}{k})\right)$ & $\text{\textnormal{cs}} \left( M_k (X_{8})\right)$ \\
\hline
3 & 0.0392668 \text{} & 0.117800 \text{} \\
 4 & 0.115898 & 0.463593 \\
 5 & 0.0209964 \text{} & 0.104982 \text{} \\
 6 & 0.149082 & 0.894495 \\
 7 & 0.0382671 \text{} & 0.267870 \text{} \\
 8 & 0.0866540 & 0.693232 \\
 9 & 0.0170042 \text{} & 0.153038 \text{} \\
 10 & 0.0636841 & 0.636841 \\
\hline
\end{tabular}
\end{tabular}

\bigskip

\begin{tabular}{ll}
\begin{tabular}{|c|c|c|}
\hline
 $k$ &  $\text{\textnormal{cs}} \left(X_{10} (\frac{2 \pi}{k})\right)$ & $\text{\textnormal{cs}} \left(M_k (X_{10})\right)$ \\
\hline
3 & 0.0749335 \text{} & 0.224800 \text{} \\
 4 & 0.218720 & 0.874878 \\
 5 & 0.0783315 \text{} & 0.391658 \text{} \\
 6 & 0.0150995 & 0.0905970 \\
 7 & 0.0560983 \text{} & 0.392688 \text{} \\
 8 & 0.0948488 & 0.758790 \\
 9 & 0.0185935 \text{} & 0.167341 \text{} \\
 10 & 0.0605490 & 0.605490 \\
\hline
\end{tabular}
&
\begin{tabular}{|c|c|c|}
\hline
 $k$ &  $\text{\textnormal{cs}} \left(X_{12} (\frac{2 \pi}{k})\right)$ & $\text{\textnormal{cs}} \left(M_k (X_{12})\right)$ \\
\hline
 3 & 0.116132 \text{} & 0.348396 \text{} \\
 4 & 0.0784470 & 0.313788 \\
 5 & 0.0428520 \text{} & 0.214260 \text{} \\
 6 & 0.0550832 & 0.330499 \\
 7 & 0.00986235 \text{} & 0.0690364 \text{} \\
 8 & 0.110442 & 0.883540 \\
 9 & 0.0276064 \text{} & 0.248458 \text{} \\
 10 & 0.0648550 & 0.648550 \\
\hline
\end{tabular}
\end{tabular}

\bigskip

\begin{tabular}{ll}
\begin{tabular}{|c|c|c|}
\hline
 $k$ &  $\text{\textnormal{cs}} \left(X_{14} (\frac{2 \pi}{k})\right)$ & $\text{\textnormal{cs}} \left(M_k (X_{14})\right)$ \\
\hline
3 & 0.161005 \text{} & 0.483014 \text{} \\
 4 & 0.192332 & 0.769328 \\
 5 & 0.0116320 \text{} & 0.0581602 \text{} \\
 6 & 0.0993703 & 0.596222 \\
 7 & 0.0393825 \text{} & 0.275677 \text{} \\
 8 & 0.00537856 & 0.0430285 \\
 9 & 0.0409719 \text{} & 0.368747 \text{} \\
 10 & 0.0735205 & 0.735205 \\
\hline
\end{tabular}
&
\begin{tabular}{|c|c|c|}
\hline
 $k$ &  $\text{\textnormal{cs}} \left(X_{16} (\frac{2 \pi}{k})\right)$ & $\text{\textnormal{cs}} \left(M_k (X_{16})\right)$ \\
\hline
 3 & 0.0416866 \text{} & 0.125060 \text{} \\
 4 & 0.0588936 & 0.235574 \\
 5 & 0.0831339 \text{} & 0.415670 \text{} \\
 6 & 0.146399 & 0.878396 \\
 7 & 0.000227239 \text{} & 0.00159067 \text{} \\
 8 & 0.0280750 & 0.224600 \\
 9 & 0.00154689 \text{} & 0.0139220 \text{} \\
 10 & 0.0849545 & 0.849545 \\
\hline
\end{tabular}
\end{tabular}

\bigskip

\begin{tabular}{|c|c|c|}
\hline
 $k$ &  $\text{\textnormal{cs}} \left(X_{18} (\frac{2 \pi}{k})\right)$ & $\text{\textnormal{cs}} \left(M_k (X_{18})\right)$ \\
\hline
  3 & 0.0907588 \text{} & 0.272277 \text{} \\
 4 & 0.177274 & 0.709096 \\
 5 & 0.0564774 \text{} & 0.282387 \text{} \\
 6 & 0.0286139 & 0.171683 \\
 7 & 0.0343586 \text{} & 0.240510 \text{} \\
 8 & 0.0526332 & 0.421066 \\
 9 & 0.0195418 \text{} & 0.175876 \text{} \\
 10 & 0.0982547 & 0.982547 \\ 
 \hline
\end{tabular}

\bigskip

\caption{Chern-Simons invariant of the hyperbolic orbifold, $\text{\textnormal{cs}} \left(X_{2n} (\frac{2 \pi}{k})\right)$ for $n$ between $1$ and $9$ and for $k$ between $3$ and $10$, and of its cyclic covering, $\text{\textnormal{cs}} \left(M_k (X_{2n})\right)$.}
\label{table1-1}
\end{table}
%%%%%%%%%%%%%%%%%%%%%%%%%%%%%%%%%%%%%%%%%%%%%%%%%%%%%%%%%
\begin{table}
\begin{tabular}{ll}
\begin{tabular}{|c|c|c|}
\hline
 $k$ &  $\text{\textnormal{cs}} \left(_{-4} (\frac{2 \pi}{k})\right)$ & $\text{\textnormal{cs}} \left(M_k (T_{-4})\right)$ \\
\hline
3 & 0.0578105 \text{} & 0.173431 \text{} \\
 4 & 0.0141698 & 0.0566791 \\
 5 & 0.0771122 \text{} & 0.385561 \text{} \\
 6 & 0.113440 & 0.680638 \\
 7 & 0.0647357 \text{} & 0.453150 \text{} \\
 8 & 0.0262590 & 0.210072 \\
 9 & 0.0506565 \text{} & 0.455908 \text{} \\
 10 & 0.0693643 & 0.693643 \\
\hline
\end{tabular}

&
\begin{tabular}{|c|c|c|}
\hline
 $k$ &  $\text{\textnormal{cs}} \left(T_{-6}(\frac{2 \pi}{k})\right)$ & $\text{\textnormal{cs}} \left(M_k (T_{-6})\right)$ \\
\hline
 3 & 0.0502767 \text{} & 0.150830 \text{} \\
 4 & 0.206063 & 0.824252 \\
 5 & 0.0724185 \text{} & 0.362092 \text{} \\
 6 & 0.136957 & 0.821740 \\
 7 & 0.0334583 \text{} & 0.234208 \text{} \\
 8 & 0.0770408 & 0.616327 \\
 9 & 0.0530941 \text{} & 0.477846 \text{} \\
 10 & 0.0324771 & 0.324771 \\
\hline
\end{tabular}
\end{tabular}

\bigskip

\begin{tabular}{ll}
\begin{tabular}{|c|c|c|}
\hline
 $k$ &  $\text{\textnormal{cs}} \left(T_{-8}(\frac{2 \pi}{k})\right)$ & $\text{\textnormal{cs}} \left(M_k (T_{-8})\right)$ \\
\hline
3 & 0.0260938 \text{} & 0.0782813 \text{} \\
 4 & 0.121024 & 0.484097 \\
 5 & 0.0343014 \text{} & 0.171507 \text{} \\
 6 & 0.123924 & 0.743545 \\
 7 & 0.0354455 \text{} & 0.248118 \text{} \\
 8 & 0.0887397 & 0.709918 \\
 9 & 0.0158804 \text{} & 0.142923 \text{} \\
 10 & 0.0555635 & 0.555635 \\
 \hline
\end{tabular}
&
\begin{tabular}{|c|c|c|}
\hline
 $k$ &  $\text{\textnormal{cs}} \left(T_{-10}(\frac{2 \pi}{k})\right)$ & $\text{\textnormal{cs}} \left(M_k (T_{-10})\right)$ \\
\hline
3 & 0.159369 \text{} & 0.478108 \text{} \\
 4 & 0.0211627 & 0.0846509 \\
 5 & 0.0799373 \text{} & 0.399686 \text{} \\
 6 & 0.0941609 & 0.564965 \\
 7 & 0.0204861 \text{} & 0.143403 \text{} \\
 8 & 0.0833782 & 0.667026 \\
 9 & 0.0170947 \text{} & 0.153852 \text{} \\
 10 & 0.0614793 & 0.614793 \\
\hline
\end{tabular}
\end{tabular}

\bigskip

\begin{tabular}{ll}
\begin{tabular}{|c|c|c|}
\hline
 $k$ &  $\text{\textnormal{cs}} \left(T_{-12}(\frac{2 \pi}{k})\right)$ & $\text{\textnormal{cs}} \left(M_k (T_{-12})\right)$ \\
\hline
3 & 0.119699 \text{} & 0.359097 \text{} \\
 4 & 0.163139 & 0.652556 \\
 5 & 0.0170874 \text{} & 0.0854371 \text{} \\
 6 & 0.0558073 & 0.334844 \\
 7 & 0.0683200 \text{} & 0.478240 \text{} \\
 8 & 0.0693583 & 0.554866 \\
 9 & 0.00963738 \text{} & 0.0867365 \text{} \\
 10 & 0.0587154 & 0.587154 \\
\hline
\end{tabular}

&
\begin{tabular}{|c|c|c|}
\hline
 $k$ &  $\text{\textnormal{cs}} \left(T_{-14}(\frac{2 \pi}{k})\right)$ & $\text{\textnormal{cs}} \left(M_k (T_{-14})\right)$ \\
\hline
 3 & 0.0758416 \text{} & 0.227525 \text{} \\
 4 & 0.0503095 & 0.201238 \\
 5 & 0.0493320 \text{} & 0.246660 \text{} \\
 6 & 0.0125167 & 0.0751005 \\
 7 & 0.0397753 \text{} & 0.278427 \text{} \\
 8 & 0.0503822 & 0.403058 \\
 9 & 0.0527761 \text{} & 0.474985 \text{} \\
 10 & 0.0509898 & 0.509898 \\
\hline
\end{tabular}
\end{tabular}

\bigskip

\begin{tabular}{ll}
\begin{tabular}{|c|c|c|}
\hline
$k$ &  $\text{\textnormal{cs}} \left(T_{-16}(\frac{2 \pi}{k})\right)$ & $\text{\textnormal{cs}} \left(M_k (T_{-16})\right)$ \\
\hline
3 & 0.0291847 \text{} & 0.0875541 \text{} \\
 4 & 0.184443 & 0.737773 \\
 5 & 0.0785018 \text{} & 0.392509 \text{} \\
 6 & 0.132806 & 0.796838 \\
 7 & 0.00813976 \text{} & 0.0569783 \text{} \\
 8 & 0.0283132 & 0.226505 \\
 9 & 0.0372653 \text{} & 0.335388 \text{} \\
 10 & 0.0401699 & 0.401699 \\
\hline
\end{tabular}
&
\begin{tabular}{|c|c|c|}
\hline
$k$ &  $\text{\textnormal{cs}} \left(T_{-18}(\frac{2 \pi}{k})\right)$ & $\text{\textnormal{cs}} \left(M_k (T_{-18})\right)$ \\
\hline
 3 & 0.147267 \text{} & 0.441800 \text{} \\
 4 & 0.0665438 & 0.266175 \\
 5 & 0.00562151 \text{} & 0.0281075 \text{} \\
 6 & 0.0843746 & 0.506248 \\
 7 & 0.0458762 \text{} & 0.321134 \text{} \\
 8 & 0.00418700 & 0.0334960 \\
 9 & 0.0196972 \text{} & 0.177275 \text{} \\
 10 & 0.0272925 & 0.272925 \\
\hline
\end{tabular}
\end{tabular}

\caption{Chern-Simons invariant of the hyperbolic orbifold, $\text{\textnormal{cs}} \left(X_{2n} (\frac{2 \pi}{k})\right)$ for $n$ between $-9$ and $-2$ and for $k$ between $3$ and $10$, and of its cyclic covering, $\text{\textnormal{cs}} \left(M_k (X_{2n})\right)$.}
\label{table1-2}
\end{table}
%%%%%%%%%%%%%%%%%%%%%%%%%%%%%%%%%%%%%%%%%%%%%%%%%%%%%%%%%%
\begin{table} 
\begin{tabular}{cc} 
\begin{tabular}{|c|c|c|}
\hline
2n & $\alpha_0$ & $\text{\textnormal{cs}}\left(X_{2n}\right)$  \\
\hline
2 & 2.40717 & 0.346796 \text{} \\
 4 & 2.75511 & 0.187220 \text{} \\
 6 & 2.87826 & 0.116482 \text{} \\
 8 & 2.94175 & 0.0787607 \text{} \\
 10 & 2.98054 & 0.0554891 \text{} \\
 12 & 3.00671 & 0.0397296 \text{} \\
 14 & 3.02556 & 0.0283589 \text{} \\
 16 & 3.03978 & 0.0197708 \text{} \\
 18 & 3.05090 & 0.0130565 \text{} \\
\hline
\end{tabular}
&
\begin{tabular}{|c|c|c|}
\hline
2n  & $\alpha_0$ & $\text{\textnormal{cs}}\left(T_{2n}\right)$  \\
\hline
 -2 & 2.09440 & 0 \\
 -4 & 2.68404 & 0.202492 \text{} \\
 -6 & 2.84713 & 0.287081 \text{} \\
 -8 & 2.92433 & 0.330333 \text{} \\
 -10 & 2.96942 & 0.356274 \text{} \\
 -12 & 2.99899 & 0.373511 \text{} \\
 -14 & 3.01989 & 0.385781 \text{} \\
 -16 & 3.03545 & 0.394957 \text{} \\
 -18 & 3.04747 & 0.402076 \text{} \\
\hline
\end{tabular}
\end{tabular}

\bigskip

\caption{Chern-Simons invariant of $X_{2n}$ for $n$ between $1$ and $9$ and for $n$ between $-9$ and $-1$).}
\label{tab2}
\end{table}
%%%%%%%%%%%%%%%%%%%%%%%%%%%%%%%%%%%%%%%%%%%%%%%%%%%%%%%%%%
%%%%%%%%%%%%%%%%%%%%%%%%%%%%%%%%%%%%%%%%%%%%%%%%%%%%%%%%%%
\medskip

\emph{Acknowledgements.} The authors would like to thank Alexander Mednykh and Hyuk Kim for their various helps.

%%%%%%%%%%%%%%%%%%%%%%%%%%%%%%%%%%%%%%%%%%%%%%%%%%%%%%%%%%%%

%%%%%%%%%%%%%%%%%%%%%%%%%%%%%%%%%%%%%%%%%%%%%%%%%%%%%%%
\end{document}